\documentclass[12pt]{article}
\usepackage{amssymb,amsmath}
\usepackage{hyperref}
\usepackage{graphicx}
\usepackage{color}

\begin{document}

\title{\LARGE\bf A rational approximation for the Dawson\text{'}s integral of real argument}

\author{
\normalsize\bf S. M. Abrarov\footnote{\scriptsize{Dept. Earth and Space Science and Engineering, York University, Toronto, Canada, M3J 1P3.}}\, and B. M. Quine$^{*}$\footnote{\scriptsize{Dept. Physics and Astronomy, York University, Toronto, Canada, M3J 1P3.}}}

\date{May, 7 2015}
\maketitle

\begin{abstract}
We present a rational approximation for the Dawson\text{'}s integral of real argument and show how it can be implemented for accurate and rapid computation of the Voigt function at small $y <  < 1$. The algorithm based on this approach enables computation with accuracy exceeding ${10^{ - 10}}$ within the domain $0 \le x \le 15$ and $0 \le y \le {10^{ - 6}}$. Due to rapid performance the proposed rational approximation runs the algorithm without deceleration.
\vspace{0.25cm}
\vspace{0.25cm}
\end{abstract}

\section{Introduction}
The Dawson\text{'}s integral is defined by \cite{Abramowitz1972, Cody1970, McCabe1974, Rybicki1989}
\begin{equation}\label{eq_1}
F\left( z \right) = {e^{ - {z^2}}}\int\limits_0^z {{e^{{t^2}}}dt},
\end{equation}
where $z = x + iy$ is a complex argument. This integral is closely related to the complex error function, also known as the Faddeeva function \cite{Abramowitz1972, Faddeyeva1961, Schreier1992}
\begin{equation}\label{eq_2}
w\left( z \right) = {e^{ - {z^2}}}\left( {1 + \frac{{2i}}{{\sqrt \pi  }}\int\limits_0^z {{e^{{t^2}}}dt} } \right).
\end{equation}
In particular, comparing equations \eqref{eq_1} and \eqref{eq_2} we can easily find a relation between the Dawson\text{'}s integral and the complex error function
\begin{equation}\label{eq_3}
F\left( z \right) = \frac{{i\sqrt \pi  }}{2}\left[ {{e^{ - {z^2}}} - w\left( z \right)} \right].
\end{equation}
As we can see from this identity, the Dawson\text{'}s integral is simply a reformulation of the complex error function.

The Dawson\text{'}s integral \eqref{eq_1} cannot be taken analytically in terms of elementary functions in a closed form and, therefore, it must be solved numerically. Several useful approximations for the Dawson\text{'}s integral have been reported  \cite{Cody1970, McCabe1974, Rybicki1989}. In this paper we show a new approach that can be successfully implemented for efficient computation.

The complex error function \eqref{eq_2} can be alternatively represented as (see Appendix A)
\begin{equation}\label{eq_4}
w\left( {x,y} \right) = \frac{1}{{\sqrt \pi  }}\int\limits_0^\infty  {\exp \left( { - {t^2}/4} \right)\exp \left( { - yt} \right)\exp \left( {ixt} \right)dt}.
\end{equation}
and expressed as a sum \cite{Schreier1992}
$$
w\left( {x,y} \right) = K\left( {x,y} \right) + iL\left( {x,y} \right),	\qquad	y \ge 0,
$$
where its real and imaginary parts are (see Appendix B)
\begin{equation}\label{eq_5}
K\left( {x,y} \right) = \frac{y}{\pi }\int\limits_{ - \infty }^\infty  {\frac{{{e^{ - {t^2}}}}}{{{y^2} + {{\left( {x - t} \right)}^2}}}dt},
\end{equation}
and
\begin{equation}\label{eq_6}
L\left( {x,y} \right) = \frac{1}{\pi }\int\limits_{ - \infty }^\infty  {\frac{{\left( {x - t} \right){e^{ - {t^2}}}}}{{{y^2} + {{\left( {x - t} \right)}^2}}}dt},
\end{equation}
respectively. The equation \eqref{eq_5} is known as the Voigt function \cite{Schreier1992, Armstrong1967, Armstrong1972, Letchworth2007, Abrarov2015a}, widely used in many branches of Applied Mathematics \cite{Srivastava1987, Abrarov2011, Pagnini2010}, Physics \cite{Edwards1992, Quine2002, Christensen2012, Berk2013, Sonnenschein2012, Quine2013} and Astronomy \cite{Emerson1996}. Since the function \eqref{eq_6} has no a specific name, it can be regarded as the $L$-function.

Due to symmetric properties of the Voigt function
$$
K\left( {x, - \left| y \right|} \right) = K\left( { - x, - \left| y \right|} \right) =  - K\left( {x,\left| y \right|} \right),
$$
it is sufficient to consider only the ${{\rm{I}}^{{\rm{st}}}}$ and ${\rm{I}}{{\rm{I}}^{{\rm{nd}}}}$ quadrants in order to cover the entire complex plane. Therefore, we will imply further that $y \ge 0$.

In our recent publication \cite{Abrarov2015a} we have shown that a new sampling methodology based on incomplete expansion of the sinc function \cite{Abrarov2015b} leads to a rational approximation of the Voigt function for rapid and accurate computation (see also the Matlab source code in \cite{Abrarov2015a})
\begin{equation}\label{eq_7}
\begin{aligned}
\kappa \left( {x,y} \right) &\buildrel \Delta \over = \sum\limits_{m = 1}^{{m_{\max }}} {\frac{{{\alpha _m}\left( {{\beta _m} + {y^2} - {x^2}} \right) + {\gamma _m}y\left( {{\beta _m} + {x^2} + {y^2}} \right)}}{{\beta _m^2 + 2{\beta _m}\left( {{y^2} - {x^2}} \right) + {{\left( {{x^2} + {y^2}} \right)}^2}}}} \\
&\Rightarrow K\left( {x,y} \right) \approx \kappa \left( {x,y + \varsigma /2} \right),
\end{aligned}
\end{equation}
where
$${\alpha _m} = \frac{{\sqrt \pi  \left( {m - 1/2} \right)}}{{2m_{\max }^2h}}\sum\limits_{n =  - N}^N {{e^{{\varsigma ^2}/4 - {n^2}{h^2}}}\sin \left( {\frac{{\pi \left( {m - 1/2} \right)\left( {nh + \varsigma /2} \right)}}{{{m_{\max }}h}}} \right)},
$$
$$
{\beta _m} = {\left( {\frac{{\pi \left( {m - 1/2} \right)}}{{2{m_{\max }}h}}} \right)^2},
$$
$$
{\gamma _m} = \frac{1}{{{m_{\max }}\sqrt \pi  }}\sum\limits_{n =  - N}^N {{e^{{\varsigma ^2}/4 - {n^2}{h^2}}}\cos \left( {\frac{{\pi \left( {m - 1/2} \right)\left( {nh + \varsigma /2} \right)}}{{{m_{\max }}h}}} \right)},
$$
$h = 0.293$, ${m_{\max }} = 12$, $\varsigma  = 2.75$ and $N = 23$ (for higher accuracy the coefficients $h$ and ${m_{\max }}$ can be taken as $0.25$ and $16$, respectively). This rational approximation is faster than the Weideman\text{'}s rational approximation \cite{Weideman1994} by factor greater than two (see \cite{Abrarov2015a} for details). Furthermore, with the same number of summation terms it is by several orders of the magnitude more accurate than the Weideman\text{'}s rational approximation in the domain of practical interest $0 < x < 40,000$ and ${10^{ - 4}} < y < {10^{ - 2}}$ \cite{Quine2002, Wells1999} required for applications using the HITRAN molecular spectroscopic database \cite{Rothman2013}. In general, the series approximation \eqref{eq_7} provides accurate results while $y \ge {10^{ - 6}}$. However, its accuracy deteriorates with decreasing parameter $y$ and, consequently, it cannot cover a narrow band region $y < {10^{ - 6}}$ along $x$-axis. It should be noted that the accuracy deterioration of the Voigt function at small $y <  < 1$ is a very common problem in most known approximations (see for example \cite{Armstrong1967, Amamou2013, Abrarov2015c}).

Although the computation of the Voigt function \eqref{eq_1} at $y < {10^{ - 6}}$ is required relatively rare in practice, it, nevertheless, has to be taken into account in algorithmic implementation. In this work we propose a rational approximation for the Dawson\text{'}s integral of real argument and show how its implementation can resolve effectively such a problem in computation of the Voigt function that occurs at small $y <  < 1$.

\section{Results and discussion}

\subsection{Derivation}

As we have shown recently, the exponential function can be expanded as a series \cite{Abrarov2015a}
\begin{equation}\label{eq_8}
\exp \left( { - {x^2}} \right) \approx \frac{1}{2}\sum\limits_{m = 1}^{{m_{\max }}} {\left[ {\frac{{{A_m} + \left( {x + i\varsigma /2} \right){B_m}}}{{C_m^2 - {{\left( {x + i\varsigma /2} \right)}^2}}} + \frac{{{A_m} + \left( { - x + i\varsigma /2} \right){B_m}}}{{C_m^2 - {{\left( { - x + i\varsigma /2} \right)}^2}}}} \right]},
\end{equation}
where
$$
{A_m} = \alpha_m = \frac{{\sqrt \pi  \left( {m - 1/2} \right)}}{{2m_{\max }^2h}}\sum\limits_{n =  - N}^N {{e^{{\varsigma ^2}/4 - {n^2}{h^2}}}\sin \left( {\frac{{\pi \left( {m - 1/2} \right)\left( {nh + \varsigma /2} \right)}}{{{m_{\max }}h}}} \right)},
$$
$$
{B_m} = - i \gamma_m = - \frac{i}{{{m_{\max }}\sqrt \pi  }}\sum\limits_{n =  - N}^N {{e^{{\varsigma ^2}/4 - {n^2}{h^2}}}\cos \left( {\frac{{\pi \left( {m - 1/2} \right)\left( {nh + \varsigma /2} \right)}}{{{m_{\max }}h}}} \right)}
$$
and
$$
{C_m^2} = \beta_m = \left( \frac{{\pi \left( {m - 1/2} \right)}}{{2{m_{\max }}h}} \right)^2.
$$
This equation has been used to obtain the rational approximation \eqref{eq_7} \cite{Abrarov2015a}. We may also utilize it for the $L$-function by substituting the exponential function approximation \eqref{eq_8} into equation \eqref{eq_6}. This leads to the following integral
\footnotesize
\begin{equation}\label{eq_9}
\begin{aligned}
L\left( {x,y}\right) &\approx \\
&\frac{1}{{2\pi }}\int\limits_{ - \infty }^\infty  {\frac{{x - t}}{{{y^2} + {{\left( {x - t} \right)}^2}}}\sum\limits_{m = 1}^{{m_{\max }}} {\left[ {\frac{{{A_m} + \left( {t + i\varsigma /2} \right){B_m}}}{{C_m^2 - {{\left( {t + i\varsigma /2} \right)}^2}}} + \frac{{{A_m} + \left( { - t + i\varsigma /2} \right){B_m}}}{{C_m^2 - {{\left( { - t + i\varsigma /2} \right)}^2}}}} \right]} } \,dt.
\end{aligned}
\end{equation}
\normalsize
The integrand of this integral is analytic everywhere over the entire complex plane except $2 + 4{m_{\max }}$ isolated points
\small
$$
\begin{aligned}
\left\{ x - iy,x + iy, - {C_m} - i\varsigma /2,{C_m} - i\varsigma /2, - {C_m} + i\varsigma /2,{C_m} \right. &+ \left. i\varsigma /2 \right\}, \\
m &\in \left\{ {1,2,3,\,\, \ldots \,\,m_{\text{max}}} \right\}
\end{aligned}
$$
\normalsize
where we can observe the singularities. In order to integrate equation (9), we may choose conveniently a contour ${C_{ccw}}$ in counterclockwise (CCW) direction on the upper-half complex plane as a semicircle centered at the origin with infinite radius. Consequently, since the domain enclosed by contour ${C_{ccw}}$ covers only the half of complex plane, the number of the isolated points is reduced twice
$$
{t_r} = \left\{ {x + iy, - {C_m} + i\varsigma /2,{C_m} + i\varsigma /2} \right\}, \quad m \in \left\{ {1,2,3,\,\, \ldots \,\,{m_{\max }}} \right\}.
$$
Using the Residue Theorem\text{'}s formula now
$$
\frac{1}{{2\pi i}}\oint\limits_{{C_{ccw}}} {f\left( t \right)} \,dt = \sum\limits_{r = 1}^{1 + 2{m_{\max }}} {{\rm{Res}}\left[ {f\left( t \right),{t_r}} \right]},
$$
where $f\left( t \right)$ is the integrand of integral \eqref{eq_9}, we obtain a rational approximation for the $L$-function \eqref{eq_6} as follows
\begin{equation}\label{eq_10}
\begin{aligned}
\lambda \left( {x,y} \right) &\buildrel \Delta \over = \sum\limits_{m = 1}^{{m_{\max }}} {\frac{{x\left[ {2{\alpha _m}y + {\gamma _m}\left( {{x^2} + {y^2} - {\beta _m}} \right)} \right]}}{{\beta _m^2 + 2{\beta _m}\left( {{y^2} - {x^2}} \right) + {{\left( {{x^2} + {y^2}} \right)}^2}}}} \\
 &\Rightarrow L\left( {x,y} \right) \approx \lambda \left( {x,y + \varsigma /2} \right).
\end{aligned}
\end{equation}
Consequently, the complex error function can be calculated as
\begin{equation}\label{eq_11}
w\left( {x,y} \right) \approx \kappa \left( {x,y + \varsigma /2} \right) + i\lambda \left( {x,y + \varsigma /2} \right).
\end{equation}

The computational testing we performed with equation \eqref{eq_11} shows that its imaginary part $\lambda \left( {x,y + \varsigma /2} \right)$ remains always accurate regardless how small the parameter $y$ is taken. Therefore the computation of the Voigt function at small $y <  < 1$ is our main objective. We may attempt to resolve this problem by using the following approximation that has been reported in our recent work \cite{Abrarov2015c}
$$
w\left( {x,y <  < 1} \right) \approx {e^{{{\left( {ix - y} \right)}^2}}}\left[ {1 + \frac{{i{e^{{x^2}}}}}{{\sqrt \pi  }}\left( {2F\left( x \right) - \frac{{1 - {e^{2ixy}}}}{x}} \right)} \right],
$$
where the Dawson\text{'}s integral of real argument is denoted as $F\left( x \right) \equiv F\left( {x,y = 0} \right)$. Taking the real part of the equation above yields the approximation of the Voigt function
\small
\begin{equation}\label{eq_12}
K\left( {x,y <  < 1} \right) \approx {e^{{y^2} - {x^2}}}\cos \left( {2xy} \right) - \frac{{2{e^{{y^2}}}}}{{\sqrt \pi  }}\left[ {y\,{\rm{sinc}}\left( {2xy} \right) - F\left( x \right)\sin \left( {2xy} \right)} \right],
\end{equation}
\normalsize
where we imply that the sinc function satisfies 
$$
\left\{ {{\rm{sinc}}\left( {2xy \ne 0} \right) = \sin \left( {2xy} \right)/\left( {2xy} \right),\,\,{\rm{sinc}}\left( {2xy = 0} \right) = 1} \right\}.
$$

According to equation \eqref{eq_3} the Dawson\text{'}s integral of complex argument can be written now as
$$
F\left( {x,y} \right) \approx \frac{{i\sqrt \pi  }}{2}\left\{ {{e^{ - {{\left( {x + iy} \right)}^2}}} - \left[ {\kappa \left( {x,y + \varsigma /2} \right) + i\lambda \left( {x,y + \varsigma /2} \right)} \right]} \right\}
$$
and since at $K\left( {x,y = 0} \right) = \exp \left( { - {t^2}} \right)$ (see Appendix B), from equation above we get an approximation for the  Dawson\text{'}s integral of real argument
\begin{equation}\label{eq_13}
F\left( {x,y = 0} \right) \approx \frac{{\sqrt \pi  }}{2}\lambda \left( {x,\varsigma /2} \right).
\end{equation}
Lastly, using the series approximation \eqref{eq_10} the equation \eqref{eq_13} can be expressed as given by
\begin{equation}\label{eq_14}
F\left( x \right) \approx \frac{{\sqrt \pi  }}{2}\sum\limits_{m = 1}^{{m_{\max }}} {\frac{{x\left[ {{\alpha _m}\varsigma  + {\gamma _m}\left( {{x^2} + {\varsigma ^2}/4 - {\beta _m}} \right)} \right]}}{{\beta _m^2 + 2{\beta _m}\left( {{\varsigma ^2}/4 - {x^2}} \right) + {{\left( {{x^2} + {\varsigma ^2}/4} \right)}^2}}}}.
\end{equation}

Figure 1 shows the difference $\varepsilon \left( t \right)$ between the original Dawson\text{'}s integral of real argument and its rational approximation \eqref{eq_14}. As we can see from this figure, despite only $12$ summation terms involved in the series approximation \eqref{eq_14}, the difference $\varepsilon \left( t \right)$ remains within a narrow range $ \pm 7 \times {10^{ - 9}}$. This confirms a rapid convergence of equation \eqref{eq_14} that makes it suitable for practical applications. Therefore, we can use it effectively as a supplement for computation of the Voigt function at small $y <  < 1$.

\begin{figure}[ht]
\begin{center}
\includegraphics[width=24pc]{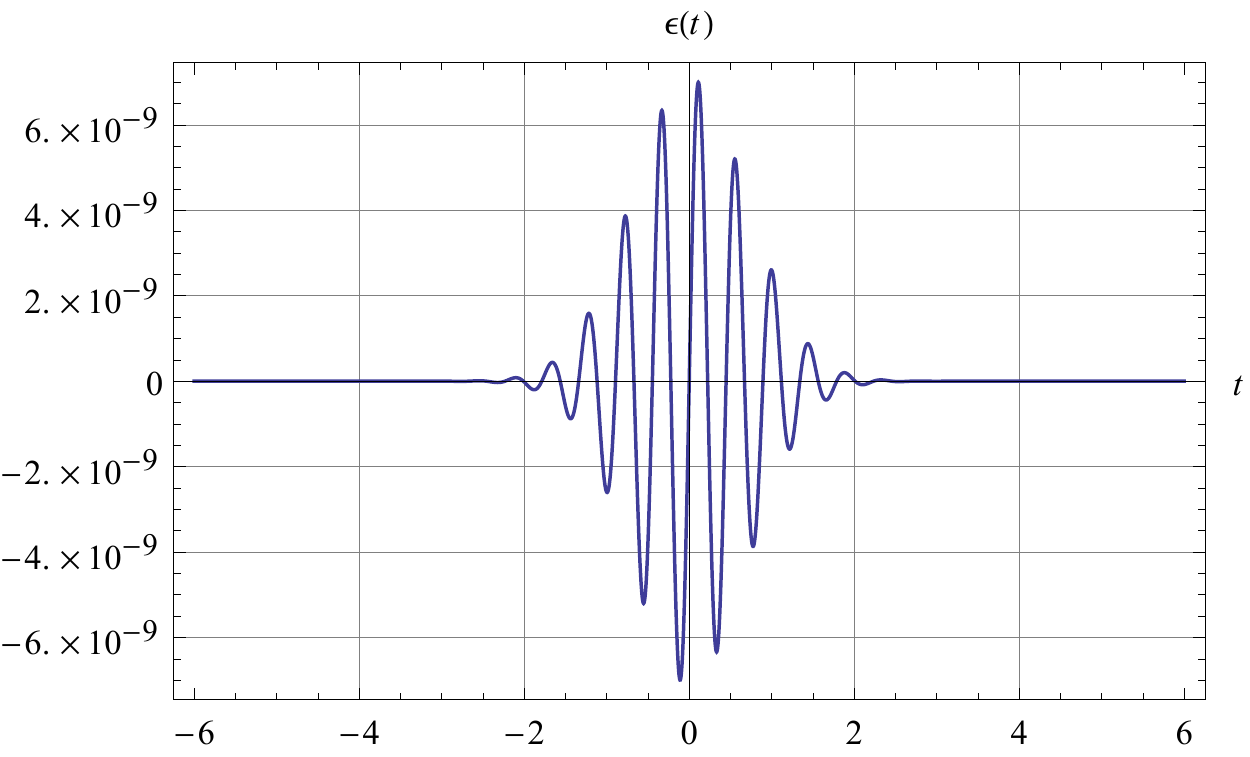}\hspace{2pc}%
\begin{minipage}[b]{28pc}
\vspace{0.3cm}
{\sffamily {\bf{Fig. 1.}} The difference $\varepsilon \left( t \right)$ at $\varsigma  = 2.75$, $N = 23$, ${m_{\max }} = 12$ and $h = 0.293$.}
\end{minipage}
\end{center}
\end{figure}

Since at $y <  < 1$ the following functions can be approximated as 
${e^{{y^2}}} \approx 1$, $\cos \left( {2xy} \right) \approx 1$, ${\rm{sinc}}\left( {2xy} \right) \approx 1$ and $\sin \left( {2xy} \right) \approx 2xy,$
the equation \eqref{eq_12} is significantly simplified as follows
$$
K\left( {x,y <  < 1} \right) \approx {e^{ - {x^2}}} - \frac{{2y}}{{\sqrt \pi  }}\left[ {1 - 2xF\left( x \right)} \right]
$$
and according to equation \eqref{eq_13} it can also be represented in form
\begin{equation}\label{eq_15}
K\left( {x,y <  < 1} \right) \approx {e^{ - {x^2}}} - \frac{{2y}}{{\sqrt \pi  }}\left[ {1 - \sqrt \pi  x\lambda \left( {x,\varsigma /2} \right)} \right].
\end{equation}

The computational time of the approximation \eqref{eq_15} is mostly taken by $\lambda $-function. Since the approximations \eqref{eq_7} and \eqref{eq_10} are about equally rapid, the approximation \eqref{eq_15} is also as fast as the approximation \eqref{eq_7}. Thus, we can see that implementation of the Dawson\text{'}s function of the real argument \eqref{eq_14} does not decelerate the computation of the Voigt function at small $y <  < 1$.

\subsection{Implementation}

When the input parameters $x$ and $y$ are large enough, say if the condition $\left| {x + iy} \right| > 15$ is satisfied, many rational approximations become effective for computation. For example, the Gauss--Hermite quadrature, the Taylor expansion series \cite{Letchworth2007} or the Laplace continued fraction \cite{Gautschi1970, Jones1988, Poppe1990a} can be used. Therefore, despite that the rational approximations \eqref{eq_7} and \eqref{eq_10} can cover accurately the entire domain of practical interest $0 < x < 40,000$ and ${10^{ - 4}} < y < {10^2}$ \cite{Quine2002, Wells1999} required for applications using the HITRAN molecular spectral database \cite{Rothman2013}, we may restrict them only within the domain $\left| {x + iy} \right| \le 15$ that is considered the most difficult for rapid and accurate computation. Thus, the computation of the Voigt function inside the domain $\left| {x + iy} \right| \le 15$ can be implemented in accordance with the following scheme
\small
$$
K\left( {x,y} \right) \approx \left\{ \begin{aligned}
&\kappa \left( {x,y + \varsigma /2} \right), \hspace{90pt} \left| {x + iy} \right| \le 15 \cap {10^{ - 6}} \le y \le 15\\
&{e^{ - {x^2}}} - \frac{{2y}}{{\sqrt \pi  }}\left[ {1 - \sqrt \pi  x\lambda \left( {x,\varsigma /2} \right)} \right], \hspace{5pt} \left| {x + iy} \right| \le 15 \cap 0 \le y < {10^{ - 6}}.
\end{aligned} \right.
$$
\normalsize

\subsection{Error analysis}

Let us define the relative error as 
\[
\Delta  = \left| {\frac{{K\left( {x,y} \right) - {K_{ref.}}\left( {x,y} \right)}}{{{K_{ref.}}\left( {x,y} \right)}}} \right|,
\]
where ${K_{ref.}}\left( {x,y} \right)$ is the reference, to quantify the accuracy of the approximation \eqref{eq_15}. The highly accurate reference values can be generated, for example, by using the Algorithm 680 \cite{Poppe1990b}, the Algorithm 916 \cite{Zaghloul2011} or a new algorithm descried in the recent work \cite{Karbach2014}.
\begin{figure}[ht]
\begin{center}
\includegraphics[width=22pc]{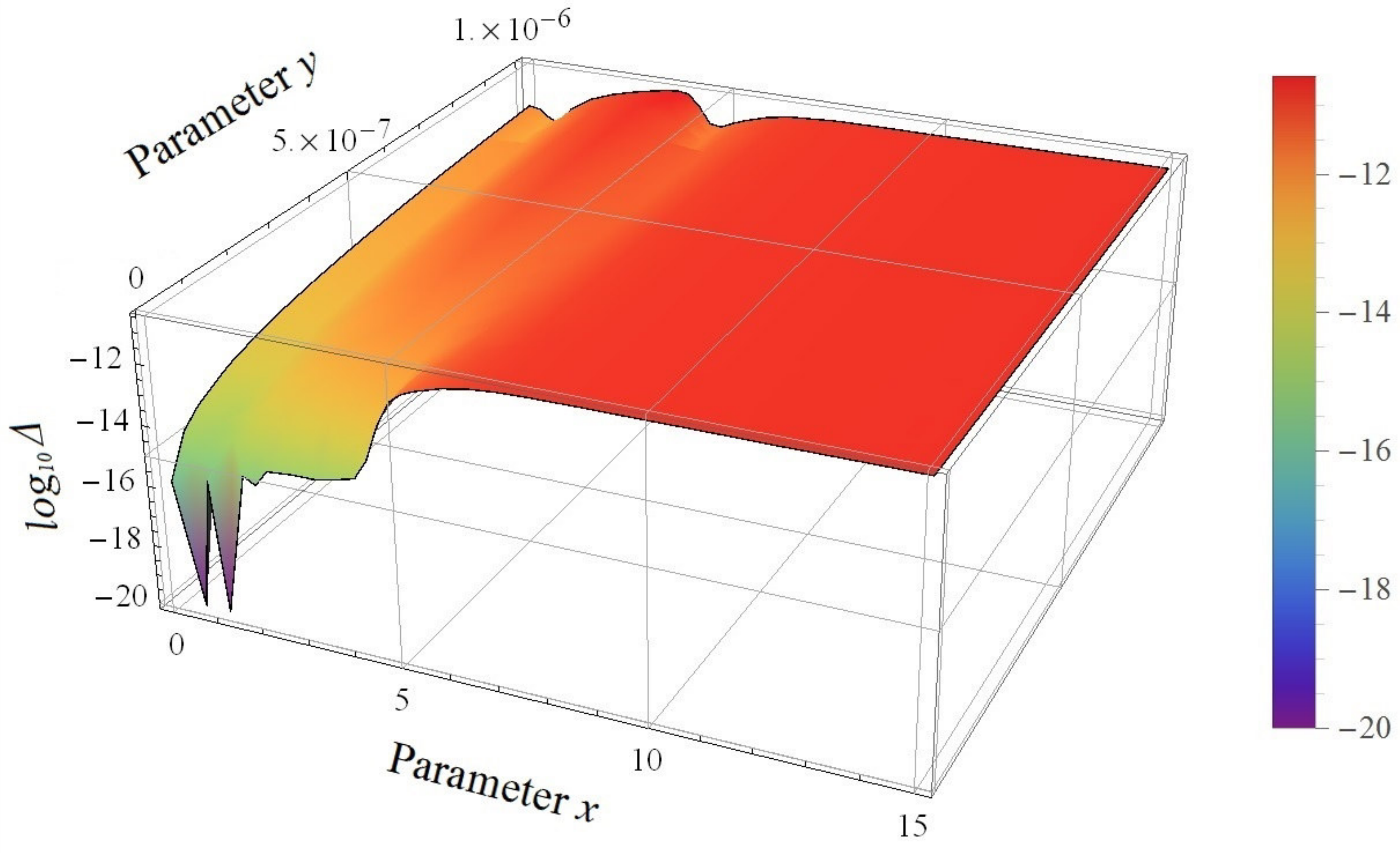}\hspace{2pc}%
\begin{minipage}[b]{28pc}
\vspace{0.3cm}
{\sffamily {\bf{Fig. 2.}} Logarithm of the relative error ${\log _{10}}\Delta $ for the Voigt function approximation \eqref{eq_15} inside the domain $0 \le x \le 15 \cap 0 \le y \le {10^{ - 6}}$ computed at $\varsigma  = 2.75$, $N = 23$, ${m_{\max }} = 12$ and $h = 0.293$.}
\end{minipage}
\end{center}
\end{figure}

Figure 2 show the logarithm of the relative error ${\log _{10}}\Delta $ for the Voigt function inside the domain $0 \le x \le 15$ and $0 \le y \le {10^{ - 6}}$. As we can see from this figure, despite that only $12$ summation terms involved in computation, the accuracy of approximation \eqref{eq_15} in this domain is better than ${10^{ - 10}}$.

Thus we can see that the approach based on implementation of the rational approximation for the Dawson\text{'}s integral of real argument \eqref{eq_14} is efficient and practical.

\section{Conclusion}

A rational approximation \eqref{eq_14} for the Dawson\text{'}s integral of real argument is presented. It can be utilized as an efficient supplement for accurate computation of the Voigt function at small $y <  < 1$. In particular, this approach enables computation with accuracy exceeding ${10^{ - 10}}$ inside the domain $0 \le x \le 15$ and $0 \le y \le {10^{ - 6}}$. Due to rapid performance the algorithmic implementation of the rational approximation \eqref{eq_14} does not decelerate the computation of the Voigt function.

\section*{Acknowledgments}

This work is supported by National Research Council Canada, Thoth Technology Inc. and York University. The authors wish to thank to Prof. Ian McDade and Dr. Brian Solheim for discussions and constructive suggestions.

\section*{Appendix A}

There is a simple proof of equation \eqref{eq_4}. The equation \eqref{eq_3} can be rearranged in form
\[
w\left( z \right) = {e^{ - {z^2}}} + \frac{{2i}}{{\sqrt \pi  }}F\left( z \right)
\]
and since
$$
{e^{ - {z^2}}} = \frac{1}{{\sqrt \pi  }}\int\limits_0^\infty  {\exp \left( { - {t^2}/4} \right)\cos \left( {zt} \right)dt},
$$
$$
F\left( z \right) = \frac{1}{2}\int\limits_0^\infty  {\exp \left( { - {t^2}/4} \right)\sin \left( {zt} \right)dt},
$$
we can write
$$
\begin{aligned}
w\left( z \right) = \frac{1}{{\sqrt \pi  }}\int\limits_0^\infty  {\exp \left( { - {t^2}/4} \right)\left( {\cos \left( {zt} \right) + i\sin \left( {zt} \right)} \right)dt} \\
 = \frac{1}{{\sqrt \pi  }}\int\limits_0^\infty  {\exp \left( { - {t^2}/4} \right) {\exp \left( {izt} \right)} dt} 
\end{aligned}
$$
or
$$
\begin{aligned}
w\left( {x,y} \right) = \frac{1}{{\sqrt \pi  }}\int\limits_0^\infty  {\exp \left( { - {t^2}/4} \right)\exp \left( {i\left( {x + iy} \right)t} \right)dt} \\
 = \frac{1}{{\sqrt \pi  }}\int\limits_0^\infty  {\exp \left( { - {t^2}/4} \right)\exp \left( { - yt} \right)\exp \left( {ixt} \right)dt} .
\end{aligned}
$$
This completes the proof.

\section*{Appendix B}

Due to discontinuity at $t = x$ the functions \eqref{eq_5} and \eqref{eq_6} cannot be defined rigorously at $y = 0$. However, since the following limits exist (see for example \cite{Armstrong1967, Letchworth2007, Zaghloul2011} for details)
$$
\mathop {\lim }\limits_{y \to 0} K\left( {x,y} \right) = {e^{ - {x^2}}},
$$
and
$$
\mathop {\lim }\limits_{y \to 0} L\left( {x,y} \right) = \frac{2}{{\sqrt \pi  }}F\left( {x,y = 0} \right)
$$
we may imply that $K\left( {x,y = 0} \right)$ and $L\left( {x,y = 0} \right)$ are defined and equal to ${e^{ - {x^2}}}$ and $\frac{2}{{\sqrt \pi  }}F\left( {x,y = 0} \right)$, respectively.



\begin{thebibliography}{9}

\bibitem{Abramowitz1972}
M. Abramowitz and I.A. Stegun. Error Function and Fresnel Integrals. Handbook of mathematical functions with formulas, graphs, and mathematical tables. $9^{\text{th}}$ ed. New York 1972, 297-309.

\bibitem{Cody1970}
W.J. Cody, K.A. Paciorek and H.C. Thacher, Chebyshev approximations for Dawson\text{'}s integral. Math. Comp. 24 (1970) 171-178. \url{http://dx.doi.org/10.1090/S0025-5718-1970-0258236-8}

\bibitem{McCabe1974}
J.H. McCabe, A continued fraction expansion with a truncation error estimate for Dawson\text{'}s integral, Math. Comp. 28 (1974) 811-816. \url{http://dx.doi.org/10.1090/S0025-5718-1974-0371020-3}

\bibitem{Rybicki1989}
G.B. Rybicki, Dawson\text{'}s integral and the sampling theorem, Comp. Phys., 3 (1989) 85-87. \url{http://dx.doi.org/10.1063/1.4822832}

\bibitem{Faddeyeva1961}
V.N. Faddeyeva, and N.M. Terent\text{'}ev, Tables of the probability integral $w\left( z \right) = {e^{ - {z^2}}}\left( {1 + \frac{{2i}}{{\sqrt \pi  }}\int_0^z {{e^{{t^2}}}dt} } \right)$ for complex argument. Pergamon Press, Oxford, 1961.

\bibitem{Schreier1992}
F. Schreier, The Voigt and complex error function: A comparison of computational methods. J. Quant. Spectrosc. Radiat. Transfer, 48 (1992) 743-762. \url{http://dx.doi.org/10.1016/0022-4073(92)90139-U}

\bibitem{Armstrong1967}
B.H. Armstrong, Spectrum line profiles: the Voigt function, J. Quant. Spectrosc. Radiat. Transfer. 7 (1967) 61-88. \url{http://dx.doi.org/10.1016/0022-4073(67)90057-X}

\bibitem{Armstrong1972}
B.H. Armstrong and B.W. Nicholls, Emission, absorption and transfer of radiation in heated atmospheres. Pergamon Press, New York, 1972.

\bibitem{Letchworth2007}
K.L. Letchworth and D.C. Benner, Rapid and accurate calculation of the Voigt function, J. Quant. Spectrosc. Radiat. Transfer, 107 (2007) 173-192. \url{http://dx.doi.org/10.1016/j.jqsrt.2007.01.052}

\bibitem{Abrarov2015a}
S.M. Abrarov and B.M. Quine, A rational approximation for efficient computation of the Voigt function in quantitative spectroscopy, arXiv:1504.00322. \url{http://arxiv.org/abs/1504.00322}

\bibitem{Srivastava1987}
H.M. Srivastava and E.A. Miller, A unified presentations of the Voigt functions, Astrophys. Space Sci., 135 (1987) 111-118. \url{http://dx.doi.org/10.1007/BF00644466}

\bibitem{Abrarov2011}
S.M. Abrarov and B.M. Quine, Efficient algorithmic implementation of the Voigt/complex error function based on exponential series approximation. Appl. Math. Comput. 218 (2011) 1894-1902. \url{http://dx.doi.org/10.1016/j.amc.2011.06.072}

\bibitem{Pagnini2010}
G. Pagnini and F. Mainardi, Evolution equations for the probabilistic generalization of the Voigt profile function, J. Comput. Appl. Math., 233 (2010) 1590-1595. \url{http://dx.doi.org/10.1016/j.cam.2008.04.040}

\bibitem{Edwards1992}
D.P. Edwards, GENLN2: A general line-by-line atmospheric transmittance and radiance model, NCAR technical note, 1992. \url{http://dx.doi.org/10.5065/D6W37T86}

\bibitem{Quine2002}
B.M. Quine and J.R. Drummond, GENSPECT: a line-by-line code with selectable interpolation error tolerance J. Quant. Spectrosc. Radiat. Transfer 74 (2002) 147-165. \url{http://dx.doi.org/10.1016/S0022-4073(01)00193-5}

\bibitem{Christensen2012}
L.E. Christensen, G.D. Spiers, R.T. Menzies and J.C Jacob, Tunable laser spectroscopy of $\rm{CO}_{2}$ near $2.05 \, {\mu}m$: Atmospheric retrieval biases due to neglecting line-mixing, J. Quant. Spectrosc. Radiat. Transfer, 113 (2012) 739-748. \url{http://dx.doi.org/10.1016/j.jqsrt.2012.02.031}

\bibitem{Berk2013}
A. Berk, Voigt equivalent widths and spectral-bin single-line transmittances: Exact expansions and the MODTRAN{\circledR}5 implementation, J. Quant. Spectrosc. Radiat. Transfer, 118 (2013) 102-120. \url{http://dx.doi.org/10.1016/j.jqsrt.2012.11.026}

\bibitem{Sonnenschein2012}
V. Sonnenschein, S. Raeder, A. Hakimi, I.D. Moore and K. Wendt, Determination of the ground-state hyperfine structure in neutral $^{229}{\rm{Th}}$, J. Phys. B: At. Mol. Opt. Phys. 45 (2012) 165005. \url{http://dx.doi.org/10.1088/0953-4075/45/16/165005}

\bibitem{Quine2013}
B.M. Quine and S.M. Abrarov, Application of the spectrally integrated Voigt function to line-by-line radiative transfer modelling. J. Quant. Spectrosc. Radiat. Transfer, 127 (2013) 37-48. \url{http://dx.doi.org/10.1016/j.jqsrt.2013.04.020}

\bibitem{Emerson1996}
D. Emerson, Interpreting astronomical spectra, John Wiley {\&} Sons Ltd, 1996.

\bibitem{Abrarov2015b}
S.M. Abrarov and B.M. Quine, Sampling by incomplete cosine expansion of the sinc function: Application to the Voigt/complex error function, Appl. Math. Comput., 258 (2015) 425-435. \url{http://dx.doi.org/10.1016/j.amc.2015.01.072}


\bibitem{Weideman1994}
J.A.C. Weideman, Computation of the complex error function. SIAM J. Numer. Anal., 31 (1994) 1497-1518. \url{http://dx.doi.org/10.1137/0731077}

\bibitem{Wells1999}
R.J. Wells, Rapid approximation to the Voigt/Faddeeva function and its derivatives. J. Quant. Spectrosc. Radiat. Transfer, 62 (1999) 29-48. \url{http://dx.doi.org/10.1016/S0022-4073(97)00231-8}

\bibitem{Rothman2013}
L.S. Rothman, I.E. Gordon, Y. Babikov, A. Barbe, D.C. Benner, P.F. Bernath, M. Birk, L. Bizzocchi, V. Boudon, L.R. Brown, A. Campargue, K. Chance, E.A. Cohen, L.H. Coudert, V.M. Devi, B.J. Drouin, A. Fayt, J.-M. Flaud, R.R. Gamache, J.J. Harrison, J.-M. Hartmann, C. Hill, J.T. Hodges, D. Jacquemart, A. Jolly, J. Lamouroux, R.J. Le Roy, G. Li, D.A. Long, O.M. Lyulin, C.J. Mackie, S.T. Massie, S. Mikhailenko, H.S.P. M\"{u}ler, O.V. Naumenko, A.V. Nikitin, J. Orphal, V. Perevalov, A. Perrin, E.R. Polovtseva and C. Richard, The HITRAN2012 molecular spectroscopic database, J. Quant. Spectrosc. Radiat. Transfer, 130 (2013) 4-50. \url{http://dx.doi.org/10.1016/j.jqsrt.2013.07.002}

\bibitem{Amamou2013}
H. Amamou, B. Ferhat and A. Bois, Calculation of the Voigt Function in the region of very small values of the parameter $a$ where the calculation is notoriously difficult, Amer. J. Anal. Chem., 4 (2013) 725-731. \url{http://dx.doi.org/10.4236/ajac.2013.412087}

\bibitem{Abrarov2015c}
S.M. Abrarov and B.M. Quine, Accurate approximations of the complex error function at small imaginary argument, J. Math. Research, 7 (2015) 44-53. \url{http://dx.doi.org/10.5539/jmr.v1n1p44}

\bibitem{Gautschi1970}
W. Gautschi, Efficient computation of the complex error function. SIAM J. Numer. Anal., 7 (1970) 187-198. \url{http://dx.doi.org/10.1137/0707012}

\bibitem{Jones1988}
W.B. Jones and W.J. Thron. Continued fractions in numerical analysis. Appl. Num. Math., 4 (1988) 143-230. http://doi.org/10.1016/0168-9274(83)90002-8

\bibitem{Poppe1990a}
G.P.M. Poppe and C.M.J. Wijers, More efficient computation of the complex error function. ACM Transact. Math. Software, 16 (1990) 38-46. \url{http://dx.doi.org/10.1145/77626.77629}

\bibitem{Poppe1990b}
G.P.M. Poppe and C.M.J. Wijers, Algorithm 680: evaluation of the complex error function. ACM Transact. Math. Software, 16 (1990) 47. \url{http://dx.doi.org/10.1145/77626.77630}

\bibitem{Zaghloul2011}
M.R. Zaghloul and A.N. Ali, Algorithm 916: computing the Faddeyeva and Voigt functions. ACM Transactions on Mathematical Software, 38 (2011) 15:1-15:22. \url{http://dx.doi.org/10.1145/2049673.2049679}

\bibitem{Karbach2014}
T.M. Karbach, G. Raven and M. Schiller, Decay time integrals in neutral meson mixing and their efficient evaluation, arXiv:1407.0748. \url{http://arxiv.org/pdf/1407.0748v1.pdf}

\end{thebibliography}
\end{document}